\documentclass[12pt,twoside]{article}
\usepackage{amssymb}

\usepackage{amsmath,amsbsy,amsfonts}
\usepackage{color}

\newtheorem{theorem}{Theorem}[section]

\let\Section=\section
\def\section{\setcounter{equation}{0}\Section}

\begin{document}

\date{}
\title{Positive solutions for infinite semipositone$\diagup $positone
quasilinear elliptic systems with singular and superlinear terms}
\author{\textsf{Brahim Khodja} \\
{\small Mathematic Department, Badji-Mokhtar Annaba University, 23000 Annaba
Algeria}\\
{\small \textit{e-mail address: }brahim.khodja@univ-annaba.org} \\
\vspace{1mm}\\
\textsf{Abdelkrim Moussaoui}\\
{\small \textit{Biology Department, A. Mira Bejaia University, }}\\
{\small \textit{Targa Ouzemour 06000 Bejaia, Algeria}}\\
{\small \textit{e-mail address:}abdelkrim.moussaoui@univ-bejaia.dz}}
\maketitle

\begin{abstract}
We establish existence and regularity of positive solutions for a class of
quasilinear elliptic systems with singular and superlinear terms. The
approach is based on sub-supersolution methods for systems of quasilinear
singular equations and the Schauder's fixed point Theorem.
\end{abstract}



\section{Introduction and main result}

\label{S1}

Let $\Omega \subset
\mathbb{R}
^{N}$ $\left( N\geq 2\right) $ is a bounded domain with $C^{1,\alpha }$%
-boundary $\partial \Omega $, $\alpha \in (0,1)$, and let $1<p,q\leq N$. We
deal with the following quasilinear singular elliptic problem%
\begin{equation}
\left\{
\begin{array}{ll}
-\Delta _{p}u=\lambda u^{\alpha _{1}}+v^{\beta _{1}} & \text{in }\Omega , \\
-\Delta _{q}v=u^{\alpha _{2}}+\lambda v^{\beta _{2}} & \text{in }\Omega , \\
u,v>0 & \text{in }\Omega , \\
u,v=0 & \text{on }\partial \Omega ,%
\end{array}%
\right.  \label{p}
\end{equation}%
where $\lambda $ is a real parameter. Here $\Delta _{p}$ and $\Delta _{q}$
denote the $p$-Laplacian and $q$-Laplacian differential operators defined by
$\Delta _{p}u=div(\left\vert \nabla u\right\vert ^{p-2}\nabla u)$ and $%
\Delta _{q}v=div(\left\vert \nabla v\right\vert ^{q-2}\nabla v)$,
respectively. We consider system (\ref{p}) in a singular case assuming that%
\begin{equation}
\begin{array}{c}
-1<\alpha _{1},\beta _{2}<0.%
\end{array}
\label{h1}
\end{equation}

We explicitly observe that under assumption (\ref{h1}) and depending on the
sign of a real number $\lambda $, it holds%
\begin{equation*}
\lim_{s\rightarrow 0^{+}}(\lambda s^{\alpha _{1}}+s^{\beta
_{1}})=\lim_{s\rightarrow 0^{+}}(s^{\alpha _{2}}+\lambda s^{\beta
_{2}})=\left\{
\begin{array}{cc}
+\infty & \text{if }\lambda >0 \\
-\infty & \text{if }\lambda <0.%
\end{array}%
\right.
\end{equation*}%
Therefore, system (\ref{p}) can be referred to as an infinite positone
problem if $\lambda >0$ and as an infinite semipositone problem if $\lambda
<0$.

The principle fact in this work is that the singularity in problem (\ref{p})
comes out through nonlinearities which are $(p-1)$-superlinear and $(q-1)$%
-superlinear near $+\infty $. Namely, we assume that%
\begin{equation}
\alpha _{2}>q-1\text{ \ and \ }\beta _{1}>p-1.  \label{h2}
\end{equation}%
In this context, system (\ref{p}) has a cooperative structure, that is, for $%
u$ (resp. $v$) fixed the right term in the first (resp. second) equation of (%
\ref{p}) is increasing in $v$ (resp. $u$). Further, according to (\ref{h2})
we have
\begin{equation*}
\lim_{s\rightarrow +\infty }(\lambda s^{\alpha _{1}}+s^{\beta
_{1}})/s^{p-1}=\lim_{s\rightarrow +\infty }(s^{\alpha _{2}}+\lambda s^{\beta
_{2}})/s^{q-1}=+\infty .
\end{equation*}

This type of problem is rare in the literature. According to our knowledge,
only a positone-type singular system with superlinear terms was examined in
\cite{YY}. There the authors considered problem (\ref{p}) depending on two
positive parameters in the whole space $%
\mathbb{R}
^{N}$. The existence of a positive entire solution is shown provided the
parameters are sufficiently small.

The sublinear condition $\alpha _{2}<q-1$ and $\beta _{1}<p-1$ for singular
systems of type (\ref{p}) have been thoroughly investigated. For a complete
overview on the study of the infinite positone problem (\ref{p}) we refer to
\cite{AC,ACG,MS,MM}, while for the study of the infinite semipositone
problem (\ref{p}), we cite \cite{EPS,LSY,LSY2}. We also mention \cite{G1,G2}
focusing on the semilinear case of (\ref{p}), that is, when $p=q=2$.

Another important class of singular problems considered in the literature is
the following%
\begin{equation}
\left\{
\begin{array}{ll}
-\Delta _{p}u=u^{\alpha _{1}}v^{\beta _{1}} & \text{in }\Omega , \\
-\Delta _{q}v=u^{\alpha _{2}}v^{\beta _{2}} & \text{in }\Omega , \\
u,v>0 & \text{in }\Omega , \\
u,v=0 & \text{on }\partial \Omega .%
\end{array}%
\right.  \label{p'}
\end{equation}%
Relevant contributions regarding the cooperative case of system (\ref{p'}),
that is $\alpha _{2},\beta _{1}>0$, can be found in \cite{GHM, GHS,MM2}.
With regard to the complementary situation $\alpha _{2},\beta _{1}<0$ which
is the so called competitive structure for system (\ref{p'}), we quote the
papers \cite{GHS,MM,MM3}. The semilinear case in (\ref{p}) (i.e. $p=q=2$)
was examined in \cite{G1,HMV,MKT} where the linearity of the principal part
is essentially used. It is worth pointing out that in the aforementioned
works, singular problem (\ref{p'}) was examined only under the sublinear
condition $\max \{\alpha _{1},\beta _{1}\}<p-1$ and $\max \{\alpha
_{2},\beta _{2}\}<q-1$. The assumptions imposed therein, especially in \cite%
{MM3}, are not satisfied for our system (\ref{p}) under hypothesis (\ref{h2}%
).

Our main concern is the question of existence of a (positive) smooth
solution for a singular system a class of elliptic systems where the
nonlinearities besides a singular terms have superlinear terms. The main
result is formulated in the next theorem.

\begin{theorem}
\label{T1}Assume (\ref{h1}) and (\ref{h2}) hold. Then system (\ref{p}) has a
(positive) solution $\left( u,v\right) $ in $C_{0}^{1,\gamma }(\overline{%
\Omega })\times C_{0}^{1,\gamma }(\overline{\Omega })$ for some $\gamma \in
(0,1)$.
\end{theorem}

The proof of Theorem \ref{T1} is done in section \ref{S2}. The main
technical difficulty consists in the presence of singular terms in system (%
\ref{p}) that can occur under hypothesis (\ref{h1}). This difficulty is
heightened by the superlinear character of (\ref{p}) that arise from (\ref%
{h2}). Our approach is chiefly based on Theorem \ref{T2} proved in Section %
\ref{S3} via Schauder's fixed point theorem (see \cite{Z}) and adequate
truncations. This is a version of the sub-supersolution method for
quasilinear singular elliptic systems with cooperative structure. We mention
that in Theorem \ref{T2} no sign condition is required on the right-hand
side nonlinearities and so it can be used for large classes of quasilinear
singular problems. A significant feature of our result lies in the obtaining
of the sub- and supersolution. Due to the superlinear character of the
nonlinearities in (\ref{p}), the latter cannot be constructed easily. At
this point, the choice of suitable functions with an adjustment of adequate
constants is crucial. Here we emphasize that the obtained sub- and
supersolution are quite different from functions considered in the
aforementioned papers, especially those constructed in \cite{MM3}.

This article is organized as follows. In section \ref{S3} we state and prove
a general theorem about sub and supersolution method for singular systems.
Section \ref{S2} contains the proof of Theorem \ref{T1}.

\section{Sub and supersolution theorem}

\label{S3}

Given $1<p<+\infty $, the spaces $L^{p}(\Omega )$ and $W_{0}^{1,p}(\Omega )$
are endowed with the usual norms $\Vert u\Vert _{p}=(\int_{\Omega }|u|^{p}\
dx)^{1/p}$ and $\Vert u\Vert _{1,p}=(\int_{\Omega }|\nabla u|^{p}\ dx)^{1/p}$%
, respectively. In the sequel, corresponding to $1<p<+\infty $, we denote $%
p^{\prime }=\frac{p-1}{p}$. We will also use the spaces $C(\overline{\Omega }%
)$ and
\begin{equation*}
C_{0}^{1,\gamma }(\overline{\Omega })=\{u\in C^{1,\gamma }(\overline{\Omega }%
):u=0\ \mbox{on
$\partial\Omega$}\}
\end{equation*}%
with $\gamma \in (0,1)$. We denote by $\lambda _{1,p}$ and $\lambda _{1,q}$
the first eigenvalue of $-\Delta _{p}$ on $W_{0}^{1,p}(\Omega )$ and of $%
-\Delta _{q}$ on $W_{0}^{1,q}(\Omega )$, respectively. Let $\phi _{1,p}$ be
the normalized positive eigenfunction of $-\Delta _{p}$ corresponding to $%
\lambda _{1,p}$, that is
\begin{equation*}
-\Delta _{p}\phi _{1,p}=\lambda _{1,p}\phi _{1,p}^{p-1}\ \text{in }\Omega ,\
\ \phi _{1,p}=0\ \text{on }\partial \Omega ,\text{ \ }\Vert \phi _{1}\Vert
_{p}=1
\end{equation*}%
Similarly, let $\phi _{1,q}$ be the normalized positive eigenfunction of $%
-\Delta _{q}$ corresponding to $\lambda _{1,q}$, that is
\begin{equation*}
-\Delta _{q}\phi _{1,q}=\lambda _{1,q}\phi _{1,q}^{q-1}\ \text{in }\Omega ,%
\text{ \ }\phi _{1,q}=0\ \text{on }\partial \Omega ,\text{ \ }\Vert \phi
_{2}\Vert _{q}=1.
\end{equation*}
For later use we set
\begin{equation}
R=\max \left\{ \underset{\overline{\Omega }}{\max }\phi _{1,p},\underset{%
\overline{\Omega }}{\max }\phi _{1,q}\right\} .  \label{40}
\end{equation}%
We denote by $d(x)$ the distance from a point $x\in \overline{\Omega }$ to
the boundary $\partial \Omega $, where $\overline{\Omega }=\Omega \cup
\partial \Omega $ is the closure of $\Omega \subset
\mathbb{R}
^{N}$. It is known that we can find a constant $l>0$ such that%
\begin{equation}
\phi _{1,p}(x),\phi _{1,q}(x)\geq ld(x)\text{ for all }x\in \Omega ,
\label{6}
\end{equation}%
where $d(x):=dist(x,\partial \Omega )$ (see, e.g., \cite{GST}).

Let us introduce the problem%
\begin{equation}
\left\{
\begin{array}{ll}
-\Delta _{p}u=f(x,u,v) & \text{in }\Omega , \\
-\Delta _{q}v=g(x,u,v) & \text{in }\Omega , \\
u,v>0 & \text{in }\Omega , \\
u,v=0 & \text{on }\partial \Omega ,%
\end{array}%
\right.  \label{p*}
\end{equation}%
where $\Omega $ is a bounded domain in $%
\mathbb{R}
^{N}$ $\left( N\geq 2\right) $ with smooth boundary, $1<p,q<\infty $ and $%
f,g:\Omega \times (0,+\infty )\times (0,+\infty )\rightarrow
\mathbb{R}
$ are continuous functions which can exhibit singularities when the
variables $u$ and $v$ approach zero. We consider system (\ref{p*}) with
cooperative structure assuming that for $u$ (resp. $v$) fixed the
nonlinearity $f$ (resp. $g$) is increasing in $v$ (resp. $u$). This makes
the sub-supersolution techniques applicable for (\ref{p*}). For systems
without cooperative structure, i.e. competitive systems, additional
assumptions are required (see \cite{GHS}).

We recall that a sub-supersolution for (\ref{p*}) is any pair $(\underline{u}%
,\underline{v})$, $(\overline{u},\overline{v})\in (W_{0}^{1,p}(\Omega )\cap
L^{\infty }(\Omega ))\times (W_{0}^{1,q}(\Omega )\cap L^{\infty }(\Omega ))$
for which there hold $(\overline{u},\overline{v})\geq (\underline{u},%
\underline{v})$ in $\Omega $,
\begin{equation*}
\begin{array}{l}
\int_{\Omega }\left\vert \nabla \underline{u}\right\vert ^{p-2}\nabla
\underline{u}\nabla \varphi \ dx-\int_{\Omega }f(x,\underline{u},\omega
_{2})\varphi \ dx \\
+\int_{\Omega }\left\vert \nabla \underline{v}\right\vert ^{q-2}\nabla
\underline{v}\nabla \psi \ dx-\int_{\Omega }g(x,\omega _{1},\underline{v}%
)\psi \ dx\leq 0,%
\end{array}%
\end{equation*}%
\begin{equation*}
\begin{array}{l}
\int_{\Omega }\left\vert \nabla \overline{u}\right\vert ^{p-2}\nabla
\overline{u}\nabla \varphi \ dx-\int_{\Omega }f(x,\overline{u},\omega
_{2})\varphi \ dx \\
+\int_{\Omega }\left\vert \nabla \overline{v}\right\vert ^{q-2}\nabla
\overline{v}\nabla \psi -\int_{\Omega }g(x,\omega _{1},\overline{v})\psi \
dx\geq 0,%
\end{array}%
\end{equation*}%
for all $\left( \varphi ,\psi \right) \in W_{0}^{1,p}\left( \Omega \right)
\times W_{0}^{1,q}\left( \Omega \right) $ with $\varphi ,\psi \geq 0$ a.e.
in $\Omega $ and for all $\left( \omega _{1},\omega _{2}\right) \in
W_{0}^{1,p}\left( \Omega \right) \times W_{0}^{1,q}\left( \Omega \right) $
satisfying $\underline{u}\leq \omega _{1}\leq \overline{u}$ and $\underline{v%
}\leq \omega _{2}\leq \overline{v}$ a.e. in $\Omega $ (see \cite[p. 269]{CLM}%
). The main goal in this section is to prove Theorem \ref{T2} below, which
is a key point in the proof of Theorem \ref{T1}.

\begin{theorem}
\label{T2}Let $\left( \underline{u},\underline{v}\right) ,$ $\left(
\overline{u},\overline{v}\right) \in C^{1}(\overline{\Omega })\times C^{1}(%
\overline{\Omega })$ be a sub and supersolution pairs of (\ref{p*}) and
suppose there exist constants $k_{1},k_{2}>0$ and $-1<\alpha ,\beta <0$ such
that%
\begin{equation}
\begin{array}{c}
\left\vert f(x,u,v)\right\vert \leq k_{1}d(x)^{\alpha }\text{ and }%
\left\vert g(x,u,v)\right\vert \leq k_{2}d(x)^{\beta }\text{ in }\Omega
\times \lbrack \underline{u},\overline{u}]\times \lbrack \underline{v},%
\overline{v}].%
\end{array}
\label{h5}
\end{equation}%
Then system (\ref{p*}) has a positive solution $(u,v)$ in $C_{0}^{1,\gamma }(%
\overline{\Omega })\times C_{0}^{1,\gamma }(\overline{\Omega })$ for certain
$\gamma \in (0,1).$
\end{theorem}

\noindent \textbf{Proof.}

For each $(z_{1},z_{2})\in C(\overline{\Omega })\times C(\overline{\Omega }),
$ let $(u,v)\in W_{0}^{1,p}(\Omega )\times W_{0}^{1,q}(\Omega )$ be the
unique solution of the problem%
\begin{equation}
\left\{
\begin{array}{ll}
-\Delta _{p}u=\widetilde{f}(x,z_{1},z_{2}) & \text{in }\Omega , \\
-\Delta _{q}v=\widetilde{g}(x,z_{1},z_{2}) & \text{in }\Omega , \\
u,v>0 & \text{in }\Omega , \\
u,v=0 & \text{on }\partial \Omega ,%
\end{array}%
\right.   \label{301}
\end{equation}%
where
\begin{equation}
\begin{array}{c}
\widetilde{f}(x,z_{1},z_{2})=f(x,\widetilde{z}_{1},\widetilde{z}_{2})\text{
and }\widetilde{g}(x,z_{1},z_{2})=g(x,\widetilde{z}_{1},\widetilde{z}_{2})%
\end{array}
\label{320}
\end{equation}%
with%
\begin{equation}
\widetilde{z}_{1}=\min (\max (z_{1},\underline{u}),\overline{u})\text{ and }%
\widetilde{z}_{2}=\min (\max (z_{2},\underline{v}),\overline{v}).
\label{300}
\end{equation}%
On account of (\ref{300}) it follows that $\underline{u}\leq \widetilde{z}%
_{1}\leq \overline{u}$ and $\underline{v}\leq \widetilde{z}_{2}\leq
\overline{v}$. Then, bearing in mind (\ref{h5}) we have%
\begin{equation}
\left\vert \widetilde{f}(x,z_{1},z_{2})\right\vert \leq k_{1}d(x)^{\alpha }%
\text{ and }\left\vert \widetilde{g}(x,z_{1},z_{2})\right\vert \leq
k_{2}d(x)^{\beta }\text{ for a.e. }x\in \Omega .  \label{302}
\end{equation}

We point out that the estimates (\ref{302}) enable us to deduce that
\begin{equation*}
\widetilde{f}(x,z_{1},z_{2})\in W^{-1,p^{\prime }}(\Omega )\text{ and }%
\widetilde{g}(x,z_{1},z_{2})\in W^{-1,q^{\prime }}(\Omega ).
\end{equation*}%
This is a consequence of Hardy-Sobolev inequality (see, e.g., \cite[Lemma 2.3%
]{AC}). Then the unique solvability of $(u,v)$ in (\ref{301}) is readily
derived from Minty-Browder Theorem (see, e.g., \cite{B}).

Let us introduce the operator
\begin{equation*}
\begin{array}{lll}
\mathcal{T}: & C(\overline{\Omega })\times C(\overline{\Omega }) &
\rightarrow C(\overline{\Omega })\times C(\overline{\Omega }) \\
& \text{ \ \ \ \ \ }(z_{1},z_{2}) & \mapsto \text{ \ \ \ }(u,v).%
\end{array}%
\end{equation*}%
We note from (\ref{301}) that the fixed point of $\mathcal{T}$ coincide with
the weak solution of (\ref{p*}). Consequently, to achieve the desired
conclusion it suffices to prove that $\mathcal{T}$ has a fixed point. To
this end we apply Schauder's fixed point theorem. Using (\ref{302}) there
exists $\gamma \in (0,1)$ such that $(u,v)\in C_{0}^{1,\gamma }(\overline{%
\Omega })\times C_{0}^{1,\gamma }(\overline{\Omega })$ and $\left\Vert
u\right\Vert _{C_{0}^{1,\gamma }(\overline{\Omega })},\left\Vert
u\right\Vert _{C_{0}^{1,\gamma }(\overline{\Omega })}\leq C,$ where $C>0$ is
independent of $u$ and $v$ (see \cite[Lemma 3.1]{H}). Then the compactness
of the embedding $C_{0}^{1,\gamma }(\overline{\Omega })\subset C(\overline{%
\Omega })$ implies that $\mathcal{T(}C(\overline{\Omega })\times C(\overline{%
\Omega }))$ is a relatively compact subset of $C(\overline{\Omega })\times C(%
\overline{\Omega })$.

Next, we show that $\mathcal{T}$ is continuous with respect to the topology
of $C(\overline{\Omega })\times C(\overline{\Omega })$. Let $%
(z_{1,n},z_{2,n})\rightarrow (z_{1},z_{2})$ in $C(\overline{\Omega })\times
C(\overline{\Omega })$ for all $n$. Denote $\left( u_{n},v_{n}\right) =%
\mathcal{T(}z_{1,n},z_{2,n})$, which reads as
\begin{equation}
\begin{array}{c}
\int_{\Omega }\left\vert \nabla u_{n}\right\vert ^{p-2}\nabla u_{n}\nabla
\varphi =\int_{\Omega }\widetilde{f}(x,z_{1,n},z_{2,n})\varphi \ dx%
\end{array}
\label{310}
\end{equation}%
and
\begin{equation}
\begin{array}{c}
\int_{\Omega }\left\vert \nabla v_{n}\right\vert ^{q-2}\nabla v_{n}\nabla
\psi =\int_{\Omega }\widetilde{g}(x,z_{1,n},z_{2,n})\psi \ dx%
\end{array}
\label{311}
\end{equation}%
for all $\left( \varphi ,\psi \right) \in W_{0}^{1,p}(\Omega )\times
W_{0}^{1,q}(\Omega )$. Inserting $(\varphi ,\psi )=(u_{n},v_{n})$ in (\ref%
{310}) and (\ref{311}), using (\ref{h5}) we get%
\begin{equation}
\left\Vert u_{n}\right\Vert _{1,p}=\int_{\Omega }\widetilde{f}%
(x,z_{1,n},z_{2,n})u_{n}\ dx\leq \int_{\Omega }d^{\alpha }u_{n}dx
\label{306}
\end{equation}%
and%
\begin{equation}
\left\Vert v_{n}\right\Vert _{1,q}=\int_{\Omega }\widetilde{g}%
(x,z_{1,n},z_{2,n})v_{n}\ dx\leq \int_{\Omega }d^{\beta }v_{n}dx.
\label{307}
\end{equation}%
Since $-1<\alpha ,\beta <0$, by virtue of the Hardy-Sobolev inequality (see,
e.g., \cite{AC}), the last integrals in (\ref{306}) and (\ref{307}) are
finite which in turn imply that $\{u_{n}\}$ and $\{v_{n}\}$ are bounded in $%
W_{0}^{1,p}(\Omega )$ and $W_{0}^{1,q}(\Omega ),$ respectively. So, passing
to relabeled subsequences, we can write the weak convergence in $%
W_{0}^{1,p}(\Omega )\times W_{0}^{1,q}\left( \Omega \right) $
\begin{equation}
\begin{array}{c}
(u_{n},v_{n})\rightharpoonup \left( u,v\right)%
\end{array}
\label{312}
\end{equation}%
for some $\left( u,v\right) \in W_{0}^{1,p}(\Omega )\times W_{0}^{1,q}\left(
\Omega \right) $. Setting $\varphi =u_{n}-u$ in (\ref{310}) and $\psi
=v_{n}-v$ in (\ref{311}), we find that
\begin{equation*}
\begin{array}{l}
\int_{\Omega }\left\vert \nabla u_{n}\right\vert ^{p-2}\nabla u_{n}\nabla
(u_{n}-u)=\int_{\Omega }\widetilde{f}(x,z_{1,n},z_{2,n})(u_{n}-u)\ dx%
\end{array}%
\end{equation*}%
and
\begin{equation*}
\begin{array}{l}
\int_{\Omega }\left\vert \nabla v_{n}\right\vert ^{p-2}\nabla v_{n}\nabla
(v_{n}-v)=\int_{\Omega }\widetilde{g}(x,z_{1,n},z_{2,n})(v_{n}-v)\ dx.%
\end{array}%
\end{equation*}%
Lebesgue's dominated convergence theorem ensures%
\begin{equation*}
\begin{array}{c}
\underset{n\rightarrow \infty }{\lim }\left\langle -\Delta
_{p}u_{n},u_{n}-u\right\rangle =\underset{n\rightarrow \infty }{\lim }%
\left\langle -\Delta _{q}v_{n},v_{n}-v\right\rangle =0.%
\end{array}%
\end{equation*}%
The $S_{+}$-property of $-\Delta _{p}$ on $W_{0}^{1,p}\left( \Omega \right) $
and of $-\Delta _{q}$ on $W_{0}^{1,q}\left( \Omega \right) $ (see, e.g. \cite%
[Proposition 3.5]{MMP}), along with (\ref{312}), implies%
\begin{equation*}
\begin{array}{c}
u_{n}\rightarrow u\text{ in }W_{0}^{1,p}(\Omega )\text{ and }%
v_{n}\rightarrow v\text{ in }W_{0}^{1,q}(\Omega ).%
\end{array}%
\end{equation*}%
Then, through (\ref{310}), (\ref{311}) and the invariance of $C(\overline{%
\Omega })\times C(\overline{\Omega })$ by $\mathcal{T}$, we infer that $%
\left( u,v\right) =\mathcal{T(}z_{1},z_{2})$. On the other hand, from (\ref%
{310}) and (\ref{311}) we know that the sequence $\{\left(
u_{n},v_{n}\right) \}$ is bounded in $C_{0}^{1,\gamma }(\overline{\Omega }%
)\times C_{0}^{1,\gamma }(\overline{\Omega })$ for certain $\gamma \in (0,1)$%
. Since the embedding $C_{0}^{1,\gamma }(\overline{\Omega })\subset C(%
\overline{\Omega })$ is compact, along a relabeled subsequence there holds $%
(u_{n},v_{n})\rightarrow (u,v)$ in $C(\overline{\Omega })\times C(\overline{%
\Omega })$. We conclude that $\mathcal{T}$ is continuous.

We are thus in a position to apply Schauder's fixed point theorem to the map
$\mathcal{T}$, which establishes the existence of $(u,v)\in C(\overline{%
\Omega })\times C(\overline{\Omega })$ satisfying $(u,v)=\mathcal{T}(u,v).$

Let us justify that%
\begin{equation*}
\underline{u}\leq u\leq \overline{u}\text{ and }\underline{v}\leq v\leq
\overline{v}\text{ in }\Omega .
\end{equation*}%
Put $\zeta =(\underline{u}-u)^{+}$ and suppose $\zeta \neq 0$. Then, bearing
in mind that system (\ref{p*}) is cooperative, from (\ref{300}), (\ref{301})
and (\ref{320}), we infer that%
\begin{equation*}
\begin{array}{c}
\int_{\{u<\underline{u}\}}|\nabla u|^{p-2}\nabla u\nabla \zeta \
dx=\int_{\Omega }|\nabla u|^{p-2}\nabla u\nabla \zeta \ dx=\int_{\{u<%
\underline{u}\}}\widetilde{f}(x,u,v)\zeta \ dx \\
\\
=\int_{\{u<\underline{u}\}}f(x,\widetilde{u},\widetilde{v})\zeta \
dx=\int_{\{u<\underline{u}\}}f(x,\underline{u},\widetilde{v})\zeta \ dx\geq
\int_{\{u<\underline{u}\}}|\nabla \underline{u}|^{p-2}\nabla \underline{u}%
\nabla \zeta \ dx.%
\end{array}%
\end{equation*}%
This implies that%
\begin{equation*}
\begin{array}{c}
\int_{\{u<\underline{u}\}}(|\nabla u|^{p-2}\nabla u-|\nabla \underline{u}%
|^{p-2}\nabla \underline{u})\nabla \zeta \ dx\leq 0,%
\end{array}%
\end{equation*}%
a contradiction. Hence $u\geq \underline{u}$ in $\Omega $. A quite similar
argument provides that $v\geq \underline{v}$ in $\Omega $. In the same way,
we prove that $u\leq \overline{u}$ and $v\leq \overline{v}$ in $\Omega $.

Finally, thanks to \cite[Lemma 3.1]{H} one has $(u,v)\in C_{0}^{1,\gamma }(%
\overline{\Omega })\times C_{0}^{1,\gamma }(\overline{\Omega })$ for some $%
\gamma \in (0,1)$. This completes the proof.
\hfill\rule{2mm}{2mm}

\section{Proof of the main result}

\label{S2}

This section is devoted to the proof of Theorem \ref{T1}. It relies on
sub-supersolution techniques shown by Theorem \ref{T2}.

Let $y_{1}$ and $y_{2}$ be the unique solutions of the problems%
\begin{equation}
\left\{
\begin{array}{ll}
-\Delta _{p}y_{1}=y_{1}^{\alpha _{1}} & \text{ in }\Omega \\
y_{1}>0 & \text{ in }\Omega \\
y_{1}=0 & \text{ on }\partial \Omega%
\end{array}%
\right. \text{ \ and \ }\left\{
\begin{array}{ll}
-\Delta _{q}y_{2}=y_{2}^{\beta _{2}} & \text{ in }\Omega \\
y_{2}>0 & \text{ in }\Omega \\
y_{2}=0 & \text{ on }\partial \Omega%
\end{array}%
\right.  \label{4}
\end{equation}%
respectively. They verify the estimates
\begin{equation}
c_{1}\phi _{1,p}(x)\leq y_{1}(x)\leq c_{2}\phi _{1,p}(x)\text{ \ and \ }%
c_{3}\phi _{1,q}(x)\leq y_{2}(x)\leq c_{4}\phi _{1,q}(x),  \label{c1}
\end{equation}%
with constants $c_{2}\geq c_{1}>0$ and $c_{4}\geq c_{3}>0$ (see \cite{GST}).
For $\delta >0$ sufficiently small we denote
\begin{equation*}
\Omega _{\delta }=\left\{ x\in \Omega :dist\left( x,\partial \Omega \right)
<\delta \right\}
\end{equation*}%
and $\mu =\mu (\delta )>0$ a constant such that%
\begin{equation}
\begin{array}{c}
\phi _{1,p}\left( x\right) ,\phi _{1,q}\left( x\right) \geq \mu \text{ in }%
\Omega \backslash \Omega _{\delta }.%
\end{array}
\label{5}
\end{equation}%
Let $\underline{u}$ and $\underline{v}$ satisfy%
\begin{equation}
-\Delta _{p}\underline{u}(x)=C\left\{
\begin{array}{ll}
y_{1}^{\alpha _{1}}(x) & \text{if \ }x\in \Omega \backslash \overline{\Omega
}_{\delta } \\
-y_{1}^{\alpha _{1}}(x) & \text{if \ }x\in \Omega _{\delta }%
\end{array}%
\right. ,\text{ }\underline{u}=0\text{ \ on }\partial \Omega  \label{3}
\end{equation}%
and%
\begin{equation}
-\Delta _{q}\underline{v}(x)=C\left\{
\begin{array}{ll}
y_{2}^{\beta _{2}}(x) & \text{if \ }x\in \Omega \backslash \overline{\Omega }%
_{\delta } \\
-y_{2}^{\beta _{2}}(x) & \text{if \ }x\in \Omega _{\delta }%
\end{array}%
\right. ,\text{ }\underline{v}=0\text{ \ on }\partial \Omega \text{,}
\label{3*}
\end{equation}%
with a constant $C>1$ to be chosen later on. The Hardy-Sobolev inequality
(see e.g. \cite{AC}) guarantees that the right hand side of (\ref{3}) and (%
\ref{3*}) are in $W^{-1,p^{\prime }}(\Omega )$ and $W^{-1,q^{\prime
}}(\Omega )$, respectively. This allows to apply the Minty-Browder theorem
(see \cite[Theorem V.15]{B}) to deduce the existence of unique solutions $%
\underline{u}$ and $\underline{v}$ for problems (\ref{3}) and (\ref{3*}),
respectively. Moreover, (\ref{4}), (\ref{c1}), (\ref{3}), (\ref{3*}) and the
monotonicity of the operators $-\Delta _{p}$ and $-\Delta _{q}$ together
with \cite[Corollary 3.1]{H} imply that%
\begin{equation}
\begin{array}{c}
\frac{c_{1}}{2}C^{\frac{1}{p-1}}\phi _{1,p}(x)\leq \underline{u}(x)\leq
c_{2}C^{\frac{1}{p-1}}\phi _{1,p}(x)\text{ and }\frac{c_{3}}{2}C^{\frac{1}{%
q-1}}\phi _{1,q}(x)\leq \underline{v}(x)\leq c_{4}C^{\frac{1}{q-1}}\phi
_{1,q}(x)\text{ in }\Omega .%
\end{array}
\label{c2}
\end{equation}%
For $\lambda \geq 0,$ the positivity of $\underline{u},\underline{v}%
,y_{1},y_{2}$ and $C$ enable us to have%
\begin{equation}
\begin{array}{c}
-Cy_{1}^{\alpha _{1}}-\lambda \underline{u}^{\alpha _{1}}\leq 0\leq
\underline{v}^{\beta _{1}}\text{ in }\Omega _{\delta }%
\end{array}
\label{14}
\end{equation}%
and%
\begin{equation}
\begin{array}{c}
-Cy_{2}^{\beta _{2}}-\lambda \underline{v}^{\beta _{2}}\leq 0\leq \underline{%
u}^{\alpha _{2}}\text{ in }\Omega _{\delta }.%
\end{array}
\label{15}
\end{equation}%
For $\lambda <0$, (\ref{c1}), (\ref{c2}) and (\ref{h1}) imply%
\begin{equation}
\begin{array}{c}
-Cy_{1}^{\alpha _{1}}-\lambda \underline{u}^{\alpha _{1}}\leq
(-Cc_{2}^{\alpha _{1}}-\lambda (\frac{c_{1}}{2}C^{\frac{1}{p-1}})^{\alpha
_{1}})\phi _{1,p}^{\alpha _{1}}\leq 0\leq \underline{v}^{\beta _{1}}\text{
in }\Omega _{\delta }%
\end{array}
\label{16}
\end{equation}%
and%
\begin{equation}
\begin{array}{c}
-Cy_{2}^{\beta _{2}}-\lambda \underline{v}^{\beta _{2}}\leq \left(
-C^{q-1}c_{4}^{\beta _{2}}-\lambda (\frac{c_{3}}{2}C^{\frac{1}{q-1}})^{\beta
_{2}}\right) \phi _{1,q}^{\beta _{2}}\leq 0\leq \underline{u}^{\alpha _{2}}%
\text{ in }\Omega _{\delta },%
\end{array}
\label{17}
\end{equation}%
provided that $C$ is sufficiently large. Now we deal with the corresponding
estimates on $\Omega \backslash \overline{\Omega }_{\delta }.$ If $\lambda
\geq 0$ we get from (\ref{c1}), (\ref{c2}), (\ref{5}), (\ref{6}) and (\ref%
{h1}) that%
\begin{equation}
\begin{array}{l}
(Cy_{1}^{\alpha _{1}}-\lambda \underline{u}^{\alpha _{1}})\underline{v}%
^{-\beta _{1}}\leq Cy_{1}^{\alpha _{1}}\underline{v}^{-\beta _{1}}\leq C^{1-%
\frac{\beta _{1}}{q-1}}c_{1}^{\alpha _{1}}(\frac{c_{3}}{2}l)^{-\beta
_{1}}\phi _{1,p}^{\alpha _{1}-\beta _{1}} \\
\leq C^{1-\frac{\beta _{1}}{q-1}}c_{1}^{\alpha _{1}}(\frac{c_{3}}{2}%
l)^{-\beta _{1}}\mu ^{\alpha _{1}-\beta _{1}}\leq 1\ \ \text{in }\Omega
\backslash \overline{\Omega }_{\delta }%
\end{array}
\label{18}
\end{equation}%
and%
\begin{equation}
\begin{array}{l}
(Cy_{2}^{\beta _{2}}-\lambda \underline{v}^{\beta _{2}})\underline{u}%
^{-\alpha _{2}}\leq Cy_{2}^{\beta _{2}}\underline{u}^{-\alpha _{2}}\leq C^{1-%
\frac{\alpha _{2}}{p-1}}c_{3}^{\beta _{2}}(\frac{c1}{2}l)^{-\alpha _{2}}\phi
_{1,q}^{\beta _{2}-\alpha _{2}} \\
\leq C^{1-\frac{\alpha _{2}}{p-1}}c_{3}^{\beta _{2}}(\frac{c1}{2}l)^{-\alpha
_{2}}\mu ^{\beta _{2}-\alpha _{2}}\leq 1\ \ \text{in }\Omega \backslash
\overline{\Omega }_{\delta },%
\end{array}
\label{19}
\end{equation}%
provided that $C$ is sufficiently large. For $\lambda <0,$ (\ref{c1}), (\ref%
{c2}), (\ref{5}), (\ref{6}) and (\ref{h1}) imply
\begin{equation*}
\begin{array}{l}
(Cy_{1}^{\alpha _{1}}-\lambda \underline{u}^{\alpha _{1}})\underline{v}%
^{-\beta _{1}}\leq \left( Cc_{1}^{\alpha _{1}}-\lambda C^{\frac{\alpha _{1}}{%
p-1}}(\frac{c_{1}}{2})^{\alpha _{1}}\right) C^{-\frac{\beta _{1}}{q-1}}(%
\frac{c_{3}}{2}l)^{-\beta _{1}}\phi _{1,p}^{\alpha _{1}-\beta _{1}} \\
=C^{1-\frac{\beta _{1}}{q-1}}c_{1}^{\alpha _{1}}\left( 1-\lambda C^{\frac{%
\alpha _{1}}{p-1}-1}2^{-\alpha _{1}}\right) (\frac{c_{3}}{2}l)^{-\beta
_{1}}\phi _{1,p}^{\alpha _{1}-\beta _{1}} \\
\leq C^{1-\frac{\beta _{1}}{q-1}}c_{1}^{\alpha _{1}}\left( 1-\lambda
2^{-\alpha _{1}}\right) (\frac{c_{3}}{2}l)^{-\beta _{1}}\mu ^{\alpha
_{1}-\beta _{1}}\leq 1\ \text{in\ }\Omega \backslash \overline{\Omega }%
_{\delta }%
\end{array}%
\end{equation*}%
and%
\begin{equation*}
\begin{array}{l}
(Cy_{2}^{\beta _{2}}-\lambda \underline{v}^{\beta _{2}})\underline{u}%
^{-\alpha _{2}}\leq \left( Cc_{3}^{\beta _{2}}-\lambda C^{\frac{\beta _{2}}{%
q-1}}(\frac{c_{3}}{2})^{\beta _{2}}\right) C^{-\frac{\alpha _{2}}{p-1}}(%
\frac{c_{1}}{2}l)^{-\alpha _{2}}\phi _{1,q}^{\beta _{2}-\alpha _{2}} \\
=C^{1-\frac{\alpha _{2}}{p-1}}c_{3}^{\beta _{2}}\left( 1-\lambda C^{\frac{%
\beta _{2}}{q-1}-1}2^{-\beta _{2}}\right) (\frac{c_{1}}{2}l)^{-\alpha
_{2}}\phi _{1,q}^{\beta _{2}-\alpha _{2}} \\
\leq C^{1-\frac{\alpha _{2}}{p-1}}c_{3}^{\beta _{2}}\left( 1-\lambda
2^{-\beta _{2}}\right) (\frac{c_{1}}{2}l)^{-\alpha _{2}}\mu ^{\beta
_{2}-\alpha _{2}}\leq 1\ \text{in\ }\Omega \backslash \overline{\Omega }%
_{\delta },%
\end{array}%
\end{equation*}%
provided that $C$ is sufficiently large. This is equivalent to%
\begin{equation}
\begin{array}{c}
Cy_{1}^{\alpha _{1}}\leq \lambda \underline{u}^{\alpha _{1}}+\underline{v}%
^{\beta _{1}}\ \ \text{in \ }\Omega \backslash \overline{\Omega }_{\delta }%
\end{array}
\label{20}
\end{equation}%
and%
\begin{equation}
\begin{array}{c}
Cy_{2}^{\beta _{2}}\leq \underline{u}^{\alpha _{2}}+\lambda \underline{v}%
^{\beta _{2}}\ \ \text{in \ }\Omega \backslash \overline{\Omega }_{\delta }%
\text{,}%
\end{array}
\label{21}
\end{equation}%
for all $\lambda \in
\mathbb{R}
.$

Due to the definition of $\underline{u}$ and $\underline{v}$ (see (\ref{3})
and (\ref{3*})) we actually have
\begin{equation}
\begin{array}{c}
\int_{\Omega }\left\vert \nabla \underline{u}\right\vert ^{p-2}\nabla
\underline{u}\nabla \varphi \text{ }dx=\int_{\Omega \backslash \overline{%
\Omega }_{\delta }}Cy_{1}^{\alpha _{1}}\varphi \text{ }dx-\int_{\Omega
_{\delta }}Cy_{1}^{\alpha _{1}}\varphi \text{ }dx%
\end{array}
\label{*}
\end{equation}%
and%
\begin{equation}
\begin{array}{c}
\int_{\Omega }\left\vert \nabla \underline{v}\right\vert ^{q-2}\nabla
\underline{v}\nabla \psi \text{ }dx=\int_{\Omega \backslash \overline{\Omega
}_{\delta }}Cy_{2}^{\beta _{2}}\psi \text{ }dx-\int_{\Omega _{\delta
}}Cy_{2}^{\beta _{2}}\psi \text{ }dx,%
\end{array}
\label{**}
\end{equation}%
where $\left( \varphi ,\psi \right) \in W_{0}^{1,p}\left( \Omega \right)
\times W_{0}^{1,q}\left( \Omega \right) $ with $\varphi ,\psi \geq 0$. Then
combining (\ref{14})-(\ref{17}), (\ref{20}), (\ref{21}) with (\ref{*}), (\ref%
{**}), it is readily seen that%
\begin{equation*}
\begin{array}{c}
\int_{\Omega }\left\vert \nabla \underline{u}\right\vert ^{p-2}\nabla
\underline{u}\nabla \varphi \text{ }dx\leq \int_{\Omega }(\lambda \underline{%
u}^{\alpha _{1}}+\underline{v}^{\beta _{1}})\varphi \text{ }dx%
\end{array}%
\end{equation*}%
and%
\begin{equation*}
\begin{array}{c}
\int_{\Omega }\left\vert \nabla \underline{v}\right\vert ^{q-2}\nabla
\underline{v}\nabla \psi \text{ }dx\leq \int_{\Omega }(\underline{u}^{\alpha
_{2}}+\lambda \underline{v}^{\beta _{2}})\psi \text{ }dx,%
\end{array}%
\end{equation*}%
for all $\left( \varphi ,\psi \right) \in W_{0}^{1,p}\left( \Omega \right)
\times W_{0}^{1,q}\left( \Omega \right) $ with $\varphi ,\psi \geq 0$,
showing that $\left( \underline{u},\underline{v}\right) $ is a subsolution
of problem (\ref{p}).

Next, we construct a supersolution part for problem (\ref{p}). To this end,
let $\widetilde{\Omega }$ be a bounded domain in $%
\mathbb{R}
^{N}$ with $C^{1,\alpha }$ boundary $\partial \widetilde{\Omega }$, $\alpha
\in (0,1)$, such that $\overline{\Omega }\subset \widetilde{\Omega }.$ We
denote by $\widetilde{\lambda }_{1,p}$ and $\widetilde{\lambda }_{1,q}$ the
first eigenvalue of $-\Delta _{p}$ on $W_{0}^{1,p}(\widetilde{\Omega })$ and
of $-\Delta _{q}$ on $W_{0}^{1,q}(\widetilde{\Omega })$, respectively. Let $%
\widetilde{\phi }_{1,p}$ be the normalized positive eigenfunction of $%
-\Delta _{p}$ corresponding to $\widetilde{\lambda }_{1,p}$, that is%
\begin{equation*}
-\Delta _{p}\widetilde{\phi }_{1,p}=\widetilde{\lambda }_{1,p}\widetilde{%
\phi }_{1,p}^{p-1}\ \text{in }\widetilde{\Omega },\ \ \widetilde{\phi }%
_{1,p}=0\text{ on }\partial \widetilde{\Omega }.
\end{equation*}%
Similarly, let $\widetilde{\phi }_{1,q}$ be the normalized positive
eigenfunction of $-\Delta _{q}$ corresponding to $\widetilde{\lambda }_{1,q}$%
, that is
\begin{equation*}
-\Delta _{q}\widetilde{\phi }_{1,q}=\widetilde{\lambda }_{1,q}\widetilde{%
\phi }_{1,q}^{q-1}\ \text{in }\Omega ,\text{ \ }\widetilde{\phi }_{1,q}=0\
\text{on }\partial \widetilde{\Omega }.
\end{equation*}%
By the definition of $\widetilde{\Omega }$ and the strong maximum principle,
there exists a constant $\rho >0$ sufficiently small such that
\begin{equation}
\widetilde{\phi }_{1,p}\left( x\right) ,\widetilde{\phi }_{1,q}\left(
x\right) >\rho \text{ in }\overline{\Omega }.  \label{10}
\end{equation}%
Without lost of generality we assume that%
\begin{equation}
R=\max \left\{ \underset{\overline{\Omega }}{\max }\widetilde{\phi }_{1,p},%
\underset{\overline{\Omega }}{\max }\widetilde{\phi }_{1,q}\right\} .
\label{9}
\end{equation}

Let $\xi _{1},\xi _{2}\in C^{1}\left( \overline{\widetilde{\Omega }}\right) $
be the solutions of the homogeneous Dirichlet problems:
\begin{equation}
\left\{
\begin{array}{l}
-\Delta _{p}\xi _{1}=C^{\delta (p-1)}\xi _{1}^{\theta _{1}}\text{ in }%
\widetilde{\Omega } \\
\xi _{1}=0\text{ on }\widetilde{\Omega }%
\end{array}%
\right. ,\text{ \ \ }\left\{
\begin{array}{l}
-\Delta _{q}\xi _{2}=C^{\delta (q-1)}\xi _{2}^{\theta _{2}}\text{ in }%
\widetilde{\Omega } \\
\xi _{2}=0\text{ on }\widetilde{\Omega },%
\end{array}%
\right.  \label{11}
\end{equation}%
with constants $\delta ,$ $\theta _{1}$ and $\theta _{2}$ satisfying%
\begin{equation}
\begin{array}{c}
\theta _{1}\in (\alpha _{1},0),\text{ \ }\theta _{2}\in (\beta _{2},0)\text{
\ and \ }\delta <\min \{\frac{1}{\theta _{1}},\frac{1}{\theta _{2}}\}<0.%
\end{array}
\label{30}
\end{equation}%
Functions $\xi _{1}$and $\xi _{2}$ verifying%
\begin{equation}
\begin{array}{c}
C^{\delta }c_{0}\widetilde{\phi }_{1,p}\leq \xi _{1}\leq C^{\delta }c%
\widetilde{\phi }_{1,p}\text{ \ and \ }C^{\delta }c_{0}^{\prime }\widetilde{%
\phi }_{1,q}\leq \xi _{2}\leq C^{\delta }c^{\prime }\widetilde{\phi }_{1,q},%
\end{array}
\label{12}
\end{equation}%
for some positive constants $c_{0},c_{0}^{\prime },c$ and $c^{\prime }$ (see
\cite{GST}). Set
\begin{equation}
\begin{array}{c}
(\overline{u},\overline{v})=C^{-\delta }(\xi _{1},\xi _{2}).%
\end{array}
\label{13}
\end{equation}%
Then we have $\left( \overline{u},\overline{v}\right) \geq \left( \underline{%
u},\underline{v}\right) $ in $\overline{\Omega }$. Indeed, on the one hand,
through (\ref{3}), (\ref{3*}) and (\ref{11}), one has
\begin{equation*}
-\Delta _{p}\overline{u}\geq -\Delta _{p}\underline{u}\text{ \ and }-\Delta
_{p}\overline{v}\geq -\Delta _{p}\underline{v}\text{ \ in }\Omega _{\delta }.
\end{equation*}%
On the other hand, on the basis of (\ref{3}), (\ref{3*}), (\ref{11}), (\ref%
{12}), (\ref{c1}), (\ref{5}), (\ref{9}), (\ref{30}) and for $C$ large
enough, we achieve%
\begin{equation*}
\begin{array}{l}
-\Delta _{p}\overline{u}=C^{-\delta (p-1)}C^{\delta (p-1)}\xi _{1}^{\theta
_{1}}=\xi _{1}^{\theta _{1}}\geq C^{\delta \theta _{1}}(c\widetilde{\phi }%
_{1,p})^{\theta _{1}} \\
\geq C^{\delta \theta _{1}}(cR)^{\theta _{1}}\geq C(c_{1}\mu )^{\alpha
_{1}}\geq C(c_{1}\phi _{1,p})^{\alpha _{1}}=Cy_{1}^{\alpha _{1}}=-\Delta _{p}%
\underline{u}\text{ \ in }\Omega \backslash \overline{\Omega }_{\delta }.%
\end{array}%
\end{equation*}%
and%
\begin{equation*}
\begin{array}{l}
-\Delta _{q}\overline{v}=C^{-\delta (q-1)}C^{\delta (q-1)}\xi _{2}^{\theta
_{2}}=\xi _{2}^{\theta _{2}}\geq C^{\delta \theta _{2}}(c^{\prime }\phi
_{1,q})^{\theta _{2}} \\
\geq C^{\delta \theta _{2}}(c^{\prime }R)^{\theta _{2}}\geq C(c_{3}\mu
)^{\beta _{2}}\geq C(c_{3}\phi _{1,q})^{\beta _{2}}\geq Cy_{2}^{\beta
_{2}}=-\Delta _{q}\underline{v}\text{ in }\Omega \backslash \overline{\Omega
}_{\delta }.%
\end{array}%
\end{equation*}%
Then the monotonicity of the operators $-\Delta _{p}$ and $-\Delta _{q}$
leads to the conclusion.

Now, taking into account (\ref{30}), (\ref{12}), (\ref{9}), (\ref{10}) and (%
\ref{h1}), for all $\lambda \in
\mathbb{R}
$, we derive that in $\overline{\Omega }$ one has%
\begin{equation}
\begin{array}{l}
\xi _{1}^{\theta _{1}}\geq C^{\delta \theta _{1}}(c\phi _{1,p})^{\theta
_{1}}\geq C^{\delta \theta _{1}}(cR)^{\theta _{1}}\geq \lambda (c_{0}\rho
)^{\alpha _{1}}+(c^{\prime }R)^{\beta _{1}} \\
\lambda (c_{0}\widetilde{\phi }_{1,p})^{\alpha _{1}}+(c^{\prime }\widetilde{%
\phi }_{1,q})^{\beta _{1}}\geq \lambda (C^{-\delta }\xi _{1})^{\alpha
_{1}}+(C^{-\delta }\xi _{2})^{\beta _{1}}=\lambda \overline{u}^{\alpha _{1}}+%
\overline{v}^{\beta _{1}}\text{ \ in }\overline{\Omega }%
\end{array}
\label{33}
\end{equation}%
and%
\begin{equation}
\begin{array}{l}
\xi _{2}^{\theta _{2}}\geq C^{\delta \theta _{2}}(c^{\prime }\phi
_{1,q})^{\theta _{2}}\geq C^{\delta \theta _{2}}(c^{\prime }R)^{\theta
_{2}}\geq \lambda (cR)^{\alpha _{2}}+(c_{0}^{\prime }\rho )^{\beta _{1}} \\
\lambda (c\widetilde{\phi }_{1,p})^{\alpha _{1}}+(c_{0}^{\prime }\widetilde{%
\phi }_{1,q})^{\beta _{1}}\geq \lambda (C^{-\delta }\xi _{1})^{\alpha
_{1}}+(C^{-\delta }\xi _{2})^{\beta _{1}}=\lambda \overline{u}^{\alpha _{1}}+%
\overline{v}^{\beta _{1}}\text{ \ in }\overline{\Omega }%
\end{array}
\label{34}
\end{equation}%
provided that $C$ is sufficiently large. Consequently, it turns out from (%
\ref{11}), (\ref{33}) and (\ref{34}) that%
\begin{equation*}
\begin{array}{l}
\int_{\Omega }\left\vert \nabla \overline{u}\right\vert ^{p-2}\nabla
\overline{u}\nabla \varphi \text{ }dx=\int_{\Omega }\xi _{1}^{\theta
_{1}}\varphi \geq \int_{\Omega }(\lambda \overline{u}^{\alpha _{1}}+%
\overline{v}^{\beta _{1}})\varphi \text{ }dx%
\end{array}%
\end{equation*}%
and%
\begin{equation*}
\begin{array}{l}
\int_{\Omega }\left\vert \nabla \overline{v}\right\vert ^{q-2}\nabla
\overline{v}\nabla \psi \text{ }dx=\int_{\Omega }\xi _{2}^{\theta _{2}}\psi
\geq \int_{\Omega }\left( \overline{u}^{\alpha _{2}}+\lambda \overline{v}%
^{\beta _{2}}\right) \psi \text{ }dx,%
\end{array}%
\end{equation*}%
for all $\left( \varphi ,\psi \right) \in W_{0}^{1,p}\left( \Omega \right)
\times W_{0}^{1,q}\left( \Omega \right) .$ This proves that the pair $(%
\overline{u},\overline{v})$ is a supersolution for problem (\ref{p}).

Finally, owing to Theorem \ref{T2} problem (\ref{p}) has a positive solution
$(u,v)\in C_{0}^{1,\gamma }(\overline{\Omega })\times C_{0}^{1,\gamma }(%
\overline{\Omega })$, for certain $\gamma \in (0,1)$, within $\left[
\underline{u},\overline{u}\right] \times \left[ \underline{v},\overline{v}%
\right] $. This completes the proof.

\end{document}